\documentstyle[12pt]{amsart}
\def\inn#1#2{\mathop \int\limits_{0}^{\infty}{#1 \rmd x}}
\newcommand{\al}{\alpha}
\newcommand{\be}{\beta}
\newcommand{\g}{\gamma}

\newcommand{\y}{\eta}

\newcommand{\lam}{\lambda}

\newcommand{\ta}{\tau}

\newcommand{\R}{{\bf R}}
\newcommand{\C}{{\bf C}}
\newcommand{\lh}{{\cal L}_2(\R)}
\newcommand{\cL}{{\cal L}_2}

\renewcommand{\Re}{{\rm Re}\;}
\renewcommand{\Im}{{\rm Im}\;}

\newcommand{\Schrodinger}{Schr\"odinger }

\newcommand{\rmd}{{\rm d}}
\newcommand{\rme}{{\rm e}}
\newcommand{\EV}{eigenvalue}
\newcommand{\EVs}{eigenvalues}
\newcommand{\EF}{eigenfunction}
\newcommand{\EFs}{eigenfunctions}
\newcommand{\sqc}{\sqrt{c}}
\newcommand{\rinf}{\rightarrow\infty}
\newcommand{\rz}{\rightarrow +0}

\begin{document}
\title[On eigenfunction approximations]{ON  EIGENFUNCTION APPROXIMATIONS\\
FOR TYPICAL NON-SELF-ADJOINT SCHR\"ODINGER OPERATORS}
\author{A. Aslanyan and E. B. Davies}
\date{17 April 1999}
\address{Department of Mathematics, King's College,
Strand, London WC2R 2LS, UK}
\email{Aslanyan@@mth.kcl.ac.uk, E.Brian.Davies@@kcl.ac.uk}
\thanks{The authors thank the Engineering
and Physical Sciences Research Council
for support under grant No. GR/L75443}
\maketitle

\begin{abstract}

We construct efficient approximations
for the \EFs \ of non-self-adjoint Schr\"o\-din\-ger operators
in one dimension. The same ideas also apply to the study of
resonances of self-adjoint \Schrodinger operators which
have dilation analytic potentials. In spite of the fact that
such eigenfunctions can have surprisingly complicated
structures with multiple local maxima, we show that a suitable adaptation of
the JWKB method is able to provide accurate global approximations
to them.
\end{abstract}


\section{Introduction}
\par
The motivation for this study was, to a great extent, our recent
paper \cite{AD}, where we discovered that extreme spectral
instability was typical for many non-self-adjoint \Schrodinger
operators, and also for complex resonances of self-adjoint
\Schrodinger operators with dilation analytic potentials.
The full analysis of such operators involves the
determination not only of their spectra, but also of their
\EFs, which may be highly non-orthogonal. It is
especially important to have a powerful tool for approximating
\EFs \ in situations in which numerical computations turn out to be
intrinsically unstable and, therefore, standard routines are
insufficiently reliable. This is the case for the spectral problems
we have been dealing with in the mentioned paper.
We start by testing the JWKB analysis of section~2 on a typical
\Schrodinger operator whose spectrum has been studied in \cite{AD}.
Then the same technique is applied to another operator
with a dilation analytic potential. Although the two operators
have different dilation analytic potentials, they exhibit surprisingly similar
spectral properties. However, in the course of this work we were primarily
interested in constructing their \EFs \ rather than comparing the corresponding
eigenvalues. In section~3 we bring together the results of computations and
report the conclusions drawn from their comparative analysis.
This allows us to suggest that the technique we have developed
can be used as a general method for a wide range of similar problems.

JWKB analysis is known to be an extremely useful technique when
applied to a range of spectral problems.
To be more precise, this notion includes a
vast variety of methods, each associated with certain asymptotic
formulae. The general asymptotic theory has been developed
in a number of books among which monographs
\cite{Fed,He,Wa} can be mentioned. There exists an extensive literature on
numerical analysis of ODEs where purely numerical procedures are combined
with and strengthened by proper JWKB methods.
Examples of different asymptotic formulae can, for instance, be found
in papers \cite{KF,KLS,Mac}.
Apart from this, one can take advantage of similar techniques for strictly
analytical purposes, including spectral analysis of \Schrodinger operators.
Recent papers on the subject include \cite{D1,D2,FK,LBM}.
In the papers \cite{D1,D2} relevant
JWKB-type formulae have been used to construct the
semi-classical modes for non-self-adjoint \Schrodinger
operators in order to study their
spectral properties. A major aspect of \cite{D1,D2}
is that the JWKB functions are defined globally, not
just asymptotically at infinity, and do not involve
analytic continuation to complex phase space.
We develop the basic approach of \cite{D2} in the next section,
working out efficient eigenfunction approximations.

Another feature of our method which makes it different from the standard
approach was suggested by our numerical results.
Having studied the spectra of several operators,
we became interested in their \EFs. The latter, although more
difficult to compute than the spectra, were also calculated
by the transfer method followed by a conventional routine.
The modes we have found are quite unexpected: some of them look
like linear combinations of several complex Gaussian
functions whose centres are quite distinct.
It is this circumstance that prompted our decision to try
to approximate the \EFs \ globally by linear combinations of several
JWKB approximate eigenfunctions. This
permitted us to approximate rather accurately
all the \EFs \ except those associated with \EVs \ which are
close to the origin,
whereas using just one JWKB function usually only works for
substantially larger \EVs. The results of section~3 show that
in our examples each \EV \ (starting from a certain number)
is associated with two JWKB functions, one of them playing
the leading role for lower and the other for higher \EVs.
Linear combinations of these
functions are proved to provide fairly good approximations
for the \EFs \ involved in our examples, in spite of the
fact that we only retain the lowest order terms in the JWKB
expansion in this paper.
Having observed this phenomenon numerically for
different examples, we propose a method suitable for
approximating \EFs \ by several easily obtained JWKB
functions. Their number can be arbitrary and depends on the
potential of the operator under consideration.

\section{JWKB Formulae}

We are interested in computing the \EVs \ and \EFs \ of the
\Schrodinger operator
\begin{equation}
H:= -\frac{\rmd^2}{\rmd x^2} + V(x)
\label{sch}
\end{equation}
where $V(x)$ is a complex-valued continuous
potential. Since we will adapt the JWKB method we study instead
the \Schrodinger operator acting in $\lh$ given by
$$
H_h:=-h^2 \frac{\rmd^2}{\rmd x^2} + V(x),
$$
where $h>0$. Throughout the paper the
potential is taken to be even; this simplification is not
necessary, but enables us to restrict our computations to
the half-line $\R_+$ instead of $\R$.

Given an \EV \  $z$ of $H_h$, let us solve the equation
\begin{equation}
(H_h - z)f(x) \ = \ 0. \label{evp}
\end{equation}
We impose one of the conditions
$$
f(0) = 0 \qquad \mbox{or} \qquad f'(0)=0
$$
bearing in mind that the \EFs \ of $H_h$ are either even
or odd. Our aim is to approximate the true
\EFs \ $f(x)$ by means of relevant JWKB asymptotic formulae.
Following \cite{D2},
we assume for now that $h$ is sufficiently small
and represent the approximate solution $y(x)$ of (\ref{evp}) as
\begin{equation}
y(x) = y(a+s) = \xi(s;h)\exp\left(-h^{-1}\psi(s)\right)
\label{ef}
\end{equation}
about some $a>0$ to be determined. One expects
$\xi(s;h)$ to be asymptotically expanded in powers of $h$:
\begin{equation}
\xi(s;h) \sim \sum_{k=0}^{\infty} h^{k}\xi_k(s), \qquad h\rz.
\label{asy}
\end{equation}
As in \cite{D2}, $\psi(s)$ is taken to be the solution of
the eikonal equation
$$
{\psi'}(s)^2 = V(a+s) - V(a) - \y^2, \qquad \psi(0) = 0,
$$
where real $\y$ and $a$ satisfy $z = \y^2 + V(a)$.
Thus, we have
\begin{equation}
\psi(s) \ = \ i\y \int_{0}^{s}\left(1 - \frac{V(a+t)-V(a)}{\y^2}\right)^{1/2}  \rmd
t.
\label{psi}
\end{equation}
Here the sign of $\y$ and the branch of the square root
in the integrand are chosen
so that $\Re\psi(s)>0$ for small $s$. An important
consequence of the fact that $V$ is complex-valued is that
generically the integrand never vanishes, and thus a unique
continuous branch of the square root is defined globally on
the real line by the above formula. Substituting
the expansion (\ref{ef}) into (\ref{evp})
and equating the coefficients by equal powers of $h$,
we get a series of equations for $\xi_k,\,k=0,1,\ldots$.
To normalise the function $y(x)$ we impose
the initial conditions $\xi_0(0)=1,\,\xi_k(0)=0,\,k>0$,
so that $y(a)=1$.
A direct calculation gives us
$$
2 \xi'_0 \psi' + \xi_0 \psi^{''} = 0
$$
and, therefore,
$$
\xi_0(s) = \left(1 - \frac{V(a+s)-V(a)}{\y^2}\right)^{-1/4}.
$$
The function $\xi_0(s)$ is globally well-defined and continuous.

The other functions $\xi_k$, $k=1,2,\ldots$, can be successively
determined as well. In this initial study, we omit the higher
order terms in (\ref{asy}) and finally obtain the approximate
solution
\begin{equation}
y(x) \ = \ \xi_0(s) \exp\left( -\frac{i\y}{h} \int_{0}^{s} \xi_0(t)^{-2}
\rmd t \right).
\label{app}
\end{equation}
We substitute $y$ directly into (\ref{evp}) to confirm that
$$
\frac{(H_h  - z)y}{y} = O(h^2), \qquad h\rz.
$$

Although the above analysis is only justified
asymptotically as $h\to +0$, we apply the formulae obtained
to the case $h=1$, that is to the operator defined by (\ref{sch}).
We refer to \cite{D2} for a
justification of this for large eigenvalues and for a wide
class of smooth potentials.  Thus we actually consider the eigenvalue problem
\begin{equation}
\left(H - \lam \right) f = 0
\label{evp1}
\end{equation}
in the rest of the paper. We expect the basic JWKB modes $y(x)$
defined by (\ref{app}) to
provide good first approximations for the \EFs \ related to
the eigenvalues $\lam$ of $H$ as $|\lam|\rinf$. The
above relation (\ref{app}) becomes the approximate
formula for the \EF \ of $H$ related to $\lam$:
\begin{equation}
y(a+s)
= \left(1 - \frac{V(a+s)-V(a)}{\y^2}\right)^{-1/4}
\exp\left( -i\y\int_{0}^{s} \left(1 - \frac{V(a+t)-V(a)}{\y^2}\right)^{1/2}
\rmd t \right)
\label{appr}
\end{equation}
where the pair of real numbers $(a,\y)$ solves
\begin{equation}
\lam = \y^2 + V(a).
\label{ei}
\end{equation}
Here $a>0$ since we consider $x \in \R_+$; the sign of $\y$
is chosen so that for small $s$ in (\ref{psi}) we have
$\Re\psi(s)>0$ for the selected branch of the square root
in the integrand. In what follows we shall use the key
formulae (\ref{appr}), (\ref{ei}) systematically. The actual efficiency
of these formulae can be judged by the
numerical results  discussed in the next section.

Note that apart from (\ref{appr}) there are other ways
of approximating $f(x)$. The expression (\ref{appr})
can, for instance,  be further approximated by
$$
\tilde{y}(x) = \left(1 - \frac{V'(a)s}{2\y^2}\right)^{-1/2}
\exp\left(-i\y s + \frac{iV'(a)s^2}{4\y}\right).
$$
This formula is easier to deal with than (\ref{appr})
because of its explicit expression. Naturally, $\tilde{y}(x)$ is expected and
proved to be less accurate for our aims than the original JWKB
function $y(x)$. One can also continue the process and work out
more terms $\xi_k,\,k=1,2,\ldots$. However, the formula
(\ref{appr}) provides quite satisfactory \EF \ approximations,
which is justified by the results of computations related to two
typical examples (see, in particular, table~6 of subsection 3.4).

Several final remarks are in order. The functions defined by
(\ref{appr}) do not generally lie in $\lh$, even when
they are good approximations to the \EFs \  in some
interval. In subsection 3.2 we discuss a further truncation
procedure needed to obtain $\lh$ functions without sacrificing
the fact that they satisfy the eigenvalue equation
approximately. Secondly the number of solutions of (\ref{ei})
for a given eigenvalue $\lam$ depends on the potential $V$,
and each of the corresponding JWKB functions may make a
contribution to the eigenfunction, as we show in subsection 3.3.
In our examples the number of solutions is always $0$, $1$ or $2$.
Finally the absence of turning points typical for the complex equation
(\ref{evp}) should be once again mentioned. It allows us
to define the functions $\psi$ and $\xi$ globally on $\R$
without applying any special complex analysis technique. This
reflects the different nature of the approach presented here
and the standard JWKB analysis.

\section{Numerical Experiments}

\subsection{JWKB Parameters}

To find the approximate \EF \ $y(x)$ we have to solve
equation (\ref{ei}) first. In this section we
consider several operators of type (\ref{sch})
starting by the
harmonic oscillator operator $H_o$ with the complex potential
$V_o(x) = (c x)^2,\,c\in\C$.
In particular,  for the $m$-th \EV \ of $H_o$ we have
the equation
$$
c (2m+1) = \y^2 + (c a)^2,  \qquad
m=0,1,\ldots ,
$$
which is solved exactly. The solution $(a,\,\y)$
such that $\y>0,\,a>0$ defines the JWKB mode $y(x)$
decaying at infinity.

As mentioned in \cite{AD}, the \EVs \
of this and similar operators are extremely unstable
under small perturbations, which makes numerical
analysis quite difficult. To be able to compute
higher \EVs \ and relevant \EFs \ numerically one
needs some additional information. For instance,
when solving \EV \  problem (\ref{evp1}) by  the
transfer method it is helpful to have at least
rough estimates for the location of ${\rm
argmax}|f(x)|$. Here we take advantage of the JWKB analysis
to find the centres of the eigenfunctions.

For comparison purposes along with equation (\ref{ei})
let us introduce its more sophisticated version:
\begin{equation}
\lam + \frac{iV'(a)}{2\y}(a,\y) - \y^2 - V(a) = 0.
\label{ti}
\end{equation}
It is obtained if we replace the \EF \
associated with an \EV \ $\lam$ by the approximation
$$
f(x) = \exp\left(-i\y s + \frac{iV'(a)s^2}{4\y}\right).
$$
For the harmonic oscillator problem we denote the roots
of (\ref{ei}) by $(a',\y')$, those of (\ref{ti})
by $(a'',\y'')$. Below we tabulate the values of
$\ta':={\rm argmax}|f(x)| - a'$ and
$\ta'':={\rm argmax}|f(x)| - a''$. The \EFs \
are computed with the use of the same basic method as
in \cite{AD}. In our numerical exercises we take a typical
value of $c=\rme^{i\pi/8}$ throughout the paper.

\pagebreak
Table~1. Parameters $a,\,\ta$ for $H_o$

\vspace{0.1in}

\begin{tabular}{|c|c|c|c|c|c|}
\hline
 $m$ & $a'$ & $a''$  & ${\rm argmax}|f_m(x)|$ & $\ta'$ & $\ta''$ \\  \hline
10& 3.371 & 3.535& 3.678   & 0.307 & 0.143\\  \hline
20& 4.711 & 4.822& 4.831   & 0.120 & 0.009\\  \hline
30& 5.746 & 5.835& 5.839   & 0.093 & 0.004\\  \hline
40& 6.621 & 6.698& 6.700   & 0.079 & 0.002\\  \hline
50& 7.393 & 7.462& 7.464   & 0.071 & 0.002\\  \hline
60& 8.092 & 8.155& 8.156   & 0.064 & 0.001\\  \hline
70& 8.736 & 8.793& 8.794   & 0.058 & 0.001\\  \hline
80& 9.335 & 9.389& 9.389   & 0.054 & 0.0\\ \hline
90& 9.897 & 9.948& 9.949   & 0.052 & 0.001\\  \hline
100& 10.430 & 10.478& 10.479  & 0.049 & 0.001\\  \hline
\end{tabular}
\vspace{0.1in}

Equation (\ref{ti}) easily solved by a standard iterative method
provides somewhat more accurate information
about the centres of the \EFs \ of $H_o$
(and other considered operators) than
(\ref{ei}) does. However, for our primary aim of
\EF \ approximation this is not vital,
and the parameters obtained from
(\ref{ei}) turn out to be more suitable. For the examples considered
below we use a standard NAG routine to find the real roots of the equations
of type (\ref{ei}) and (\ref{ti}).

An interesting example where we apply the JWKB method
is an operator with a dilation analytic potential
studied in \cite{AD}. This operator is given by
$$
H_c = - c^{-1}\frac{\rmd^2}{\rmd x^2} + V(\sqc x),
\qquad V(x)= x^2\exp(-x^2/b^2),
$$
where $c\in\C$, $b\in\R$ as in \cite{AD}.
Equation (\ref{ei}) then becomes
\begin{equation}
\lam_m c = \y^2 + c V(\sqc a)
\label{eik}
\end{equation}
where $\lam_m$ denotes the $m$-th eigenvalue.
Another typical quantum mechanical operator $\tilde{H}_c$
to be considered here has a similar form to $H_c$:
$$
\tilde{H}_c = - c^{-1}\frac{\rmd^2}{\rmd x^2} + \tilde{V}(\sqc x),
$$
$$
\tilde{V}(x) = \al \left(\exp(-\g(x-\be)^2) + \exp(-\g(x+\be)^2)\right),
 \qquad \al, \be, \g \in \R.
$$
Both $H_c$ and $\tilde{H}_c$ are the operator families
parametrised by a complex $c$; their \EVs \ are known to
be independent of $c$ in a sense explained in \cite{AD}.
If one puts $c=1$ then the
corresponding operators $H_1$ and $\tilde{H}_1$
are self-adjoint and known to have
complex resonances which are the \EVs \ of the original
non-self-adjoint operators (see, for example, \cite{CFKS} where
the theory of dilation analytic resonances is exposed).

The spectra of  $H_o$ and $H_c$ for
sufficiently large $b$ are proved to be very close to each other;
the same can be said about the \EFs \ of the two operators. For $b=100$
let us tabulate the JWKB parameters $a''$ found from
(\ref{ti}) (where obvious changes are made to suit $H_c$)
to make sure they also approximate the actual values of ${\rm argmax}|f_m(x)|$.

\vspace{0.1in}

Table~2. Parameters $a'',\,\ta''$ for $H_c$; $b=100$

\vspace{0.1in}

\begin{tabular}{|c|c|c|c|c|c|}
\hline
 $m$ & $a''$  & ${\rm argmax}|f_m(x)|$ & $\ta''$ \\  \hline
10& 3.537& 3.682   & 0.145 \\  \hline
20& 4.826& 4.836   & 0.010\\  \hline
30& 5.842& 5.850   & 0.008\\  \hline
40& 6.708& 6.713   & 0.005\\  \hline
50& 7.475& 7.479   & 0.004\\  \hline
60& 8.172& 8.175   & 0.003\\  \hline
70& 8.815& 8.817   & 0.002\\  \hline
80& 9.414& 9.415   & 0.001\\  \hline
90& 9.982& 9.983   & 0.001\\  \hline
100& 10.519& 10.520 & 0.001\\  \hline
\end{tabular}

\vspace{0.1in}

Comparing tables 1 and 2 one can see that
the results obtained for $H_o$ and $H_c$ ($b=100$)
are quite similar. The plots of the computed
\EFs \ indicate the same fact. The comparison of the actual
\EFs \ of these two operators contributes to our study of
their spectra. From the numerical point of view it may turn out
to be useful to have the values of the parameter $a''$ found
from (\ref{ti}) which approximate the centres
of the eigenfunctions reasonably accurately.

\begin{figure}[h]
\begin{picture}(100,100)(0,0)
\includegraphics{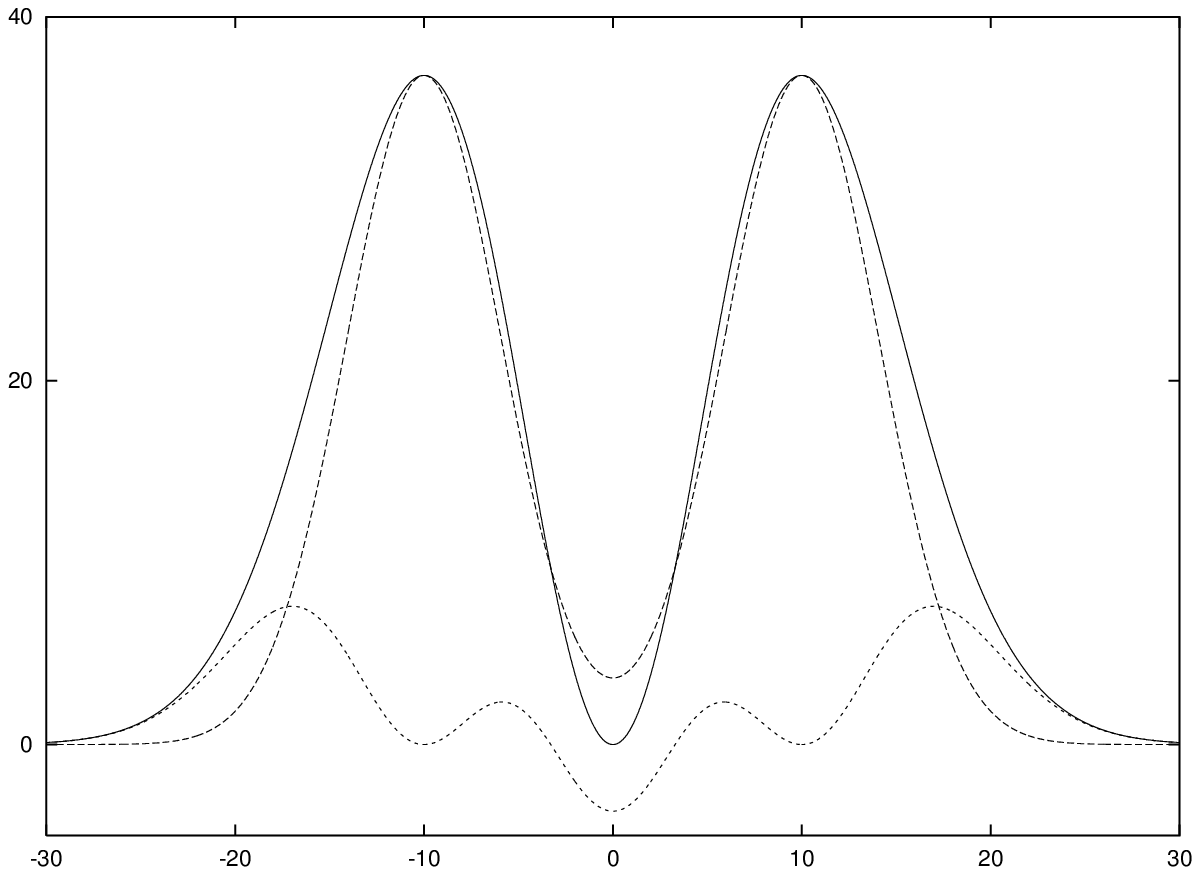}
\end{picture}
\end{figure}
\vspace{.3in}

\begin{center}
Figure 1. Potentials $V(x)$ (solid line) and
$\tilde{V}(x)$ (dotted line); values of $V-\tilde{V}$
\end{center}

From now on we shall concentrate on the two operators
$H_c$ at $b=10$ and $\tilde{H}_c$ where $\g=0.03$, $\al=100/\rme$, $\be=10$.
The choice of the parameters $\al,\,\be,\,\g$
involved into $\tilde{H}_c$ is systematic rather than random.
Our idea is to take two close self-adjoint operators $H_1$ and $\tilde{H}_1$
and transform them into non-self-adjoint operators by the
dilation analyticity technique (as has been done
for $H_1$ in \cite{AD}). In fact, the potentials
$V(x)$ and $\tilde{V}(x)$ at the chosen values of
$b$ and $\al,\,\be,\,\g$ are relatively close to one another:
$\| V(x)-\tilde{V}(x)\|/\|V(x)\| \approx 0.18$; see also figure~1
where $V$, $\tilde{V}$ and $V-\tilde{V}$ are plotted.
(Here and below by $\|\cdot\|$ we understand the
$\cL$-norm.) Of course, the
complex potentials $V(\sqc x)$ and $\tilde{V}(\sqc x)$
thus obtained are no longer similar to each other.
Remarkably, the operators $H_c$ and $\tilde{H}_c$, although
different, prove to have similar spectral properties.
Namely, their lower \EVs \ are real (up to a chosen accuracy);
starting from a certain number they turn sharply into the
lower half-plane. This behaviour of the complex resonances
observed in \cite{AD,Br} seems to be typical for a
range of operators, and the results
quoted here provide yet another numerical evidence of
the fact.

Below we quote the even \EVs \ of $\tilde{H}_c$
calculated by means of the same procedure as proposed and
implemented in our earlier paper. The \EVs \ of $H_c$
are tabulated in \cite{AD} (see table~10 there and also figure~2 below for comparisons).

\vspace{.1in}

Table~3. Eigenvalues $\tilde{\lam}_m$ of $\tilde{H}_c$

\vspace{.1in}

\begin{tabular}{|c|c|c|c|c|c|c|c|}
\hline
$m$ & $\tilde{\lam}_m$& \ & $m$ & $\tilde{\lam}_m$& \ & $m$ & $\tilde{\lam}_m$\\ \hline
0  &  4.4063         & \ & 18 & 30.2243 & \ &36 &  $46.9295  -7.3270i$\\ \hline
2  &  7.3902         & \ & 20 & 32.6595 & \ & 38 &  $48.7860  -9.0639i$\\ \hline
4  & 10.3846         & \ & 22 & $34.8766  -0.0005i$& \ & 40 &  $50.6604 -10.9142i$\\ \hline
6  & 13.3718         & \ & 24 & $36.6986  -0.0736i$& \ & 42 &  $52.5485 -12.8733i$\\ \hline
8  & 16.3358         & \ & 26 & $38.1852  -0.6680i$& \ & 44 &  $54.4504 -14.9375i$\\ \hline
10 & 19.2609         & \ & 28 &  $39.8085  -1.6721i$& \ & 46 &  $56.3613 -17.1037i$\\ \hline
12 & 22.1312         & \ & 30 &  $41.5232  -2.8670i$& \ & 48 &  $58.2803 -19.3695i$\\ \hline
14 & 24.9295         & \ & 32 &  $43.2907  -4.2188i$& \ & 50 &  $60.2058 -21.7326i$\\ \hline
16 & 27.6357         & \ & 34 &  $45.0959  -5.7093i$& \ & 52 &  $62.1367 -24.1910i$\\ \hline
\end{tabular}

\vspace{.1in}

Given the \EVs, we find the values of $a$ and $\y$
from (\ref{eik}). In a generic situation this equation
may have several real
solutions $(a_k,\,\y_k)$ such that $a_k>0$ and $\y$ has
the appropriate sign as discussed in section~2. This is
the case for both of the operators $H_c$ and $\tilde{H}_c$,
unlike the former two operators. Namely, for higher \EVs \  there exist two
solutions providing $|y(a_k+s)|<|y(a_k)|$ for sufficiently small
$s$. The values of $a_k,\,\y_k,\,k=1,2,\,a_2>a_1$, are tabulated below (see table~4).
The results for $H_c$ and $\tilde{H}_c$
are qualitatively similar; the second solution of (\ref{eik})
appears for the 28-th and 30-th eigenvalue respectively.

\vspace{.1in}

Table~4. Solutions of (\ref{eik}): $a,\,\y$ relate to $H_c$,
$\tilde{a},\,\tilde{\y}$ to $\tilde{H}_c$

\vspace{.1in}
\begin{tabular}{|c|c|c|c|c|c|c|c|c|c|c|} \hline
  $m$ & $a$ & $\y$ & $\tilde{a}$ & $\tilde{\y}$
& \ & $m$  & $a$ & $\y$ & $\tilde{a}$ & $\tilde{\y}$\\ \hline
20  & 5.0880 &3.9262 &  5.4144  & 3.8701 & &
40  & 3.5726 & 6.3897  &  4.4905 & 6.2876 \\
26  & 5.4023 & 4.3319 &  5.8455  & 4.2843 &
 &   & 12.3112 & -3.5569  & 11.4371 &-3.3463 \\ \hline
28  & 5.2223 & 4.5915 &  5.7562  & 4.5455 & &
42  & 3.2224 & 6.6574   &  4.1567 & 6.5730 \\
 & 10.7554 & -0.3678  & & &
& & 12.6198 & -3.9050    & 11.6097 &-3.7166 \\ \hline
30  & 5.0489 & 4.8785  &  5.6270 & 4.8273 & &
44  & 2.7926  & 6.9310  &  3.7719 & 6.8549 \\
 & 10.9252 & -1.3639   & 10.7532 &-0.1320 &
& & 12.9812 & -4.2601    & 11.7924 &-4.0723 \\ \hline
32  & 4.7845 & 5.1944  &  5.4631 & 5.1183 & &
46  &  2.2971 & 7.1949  &  3.3172 & 7.1332 \\
    & 11.1804 & -1.9060  & 10.8609 &-1.4247 &
 &   & 13.3664 & -4.6127   & 11.9855 &-4.4182 \\ \hline
34  &   4.5351 & 5.5016 &  5.2676 & 5.4124 & &
48  & 1.6747 & 7.4508   &  2.7555 & 7.4085 \\
 & 11.4189 & -2.3884    & 10.9854 &-2.0364 &
 & & 13.7846 & -4.9684   & 12.1898 &-4.7583\\ \hline
36  & 4.2522 & 5.8023  &  5.0413 & 5.7065 & &
50 & 0.6079 & 7.7002 & 1.9911 & 7.6815 \\
    & 11.6864 & -2.8042  & 11.1241 &-2.5247 &
 & &  14.2496 & -5.3334  & 12.4065 & -5.0958  \\ \hline
38  & 3.9280 & 6.0983  &  4.7833 & 5.9985 & &
52 &  &                & 0.2650 & 7.9532 \\
 &  11.9881 & -3.1856  & 11.2750 &-2.9532 &
 & &  &  &  12.6375 & -5.4334 \\ \hline
\end{tabular}
\vspace{.1in}

Figure 2 explains why the second solution $(a,\eta)$ appears. Here
the complex potentials $cV(\sqc x)$ and $c\tilde{V}(\sqc x)$ are plotted for
$0 \leq x \leq X$, $X$ being sufficiently large. We treat $x$ as a
parameter and represent the data in the form $(\Re(cV),\,\Im(cV))$ for each $x$.
The dots in the figure denote the values of $\lam_mc$, i.e.,
the \EVs \ of the operators $H_c$ and $\tilde{H}_c$ multiplied
by $c$ (only even \EVs \ are pictured for both operators). These pictures refer to
equation (\ref{eik}) and show how it is solved graphically.
Clearly, if a dot $\lam_mc$ is `inside' the graph then the two possible values of
$\y^2=\Re(\lam_mc)-\Re(cV)$ have opposite signs, while dots lying
`outside' determine two positive values of $\y^2$.
Hence in the former case there is only one real
pair $(a,\y)$ solving (\ref{eik}) (provided the proper sign of $\y$ is chosen),
whereas in the latter (\ref{eik}) has two solutions of this kind.
As is seen, the critical points are $m=28$ for the first operator
and $m=30$ for the second. It is from these numbers onwards that the
behaviour of \EFs \ changes (see figure~3 where
\EFs \ of $H_c$ are plotted).

\begin{figure}[h]
\begin{picture}(200,170)(0,0)
\includegraphics{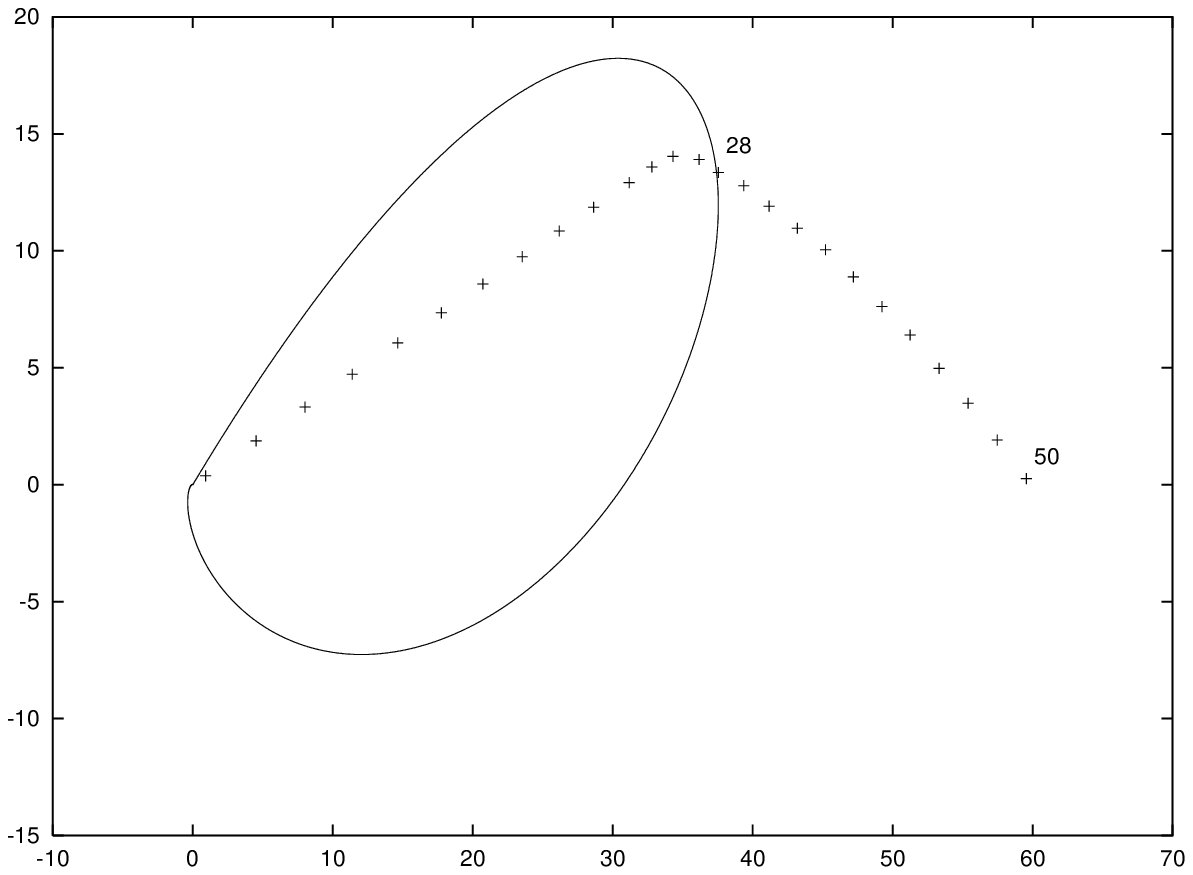}
\includegraphics{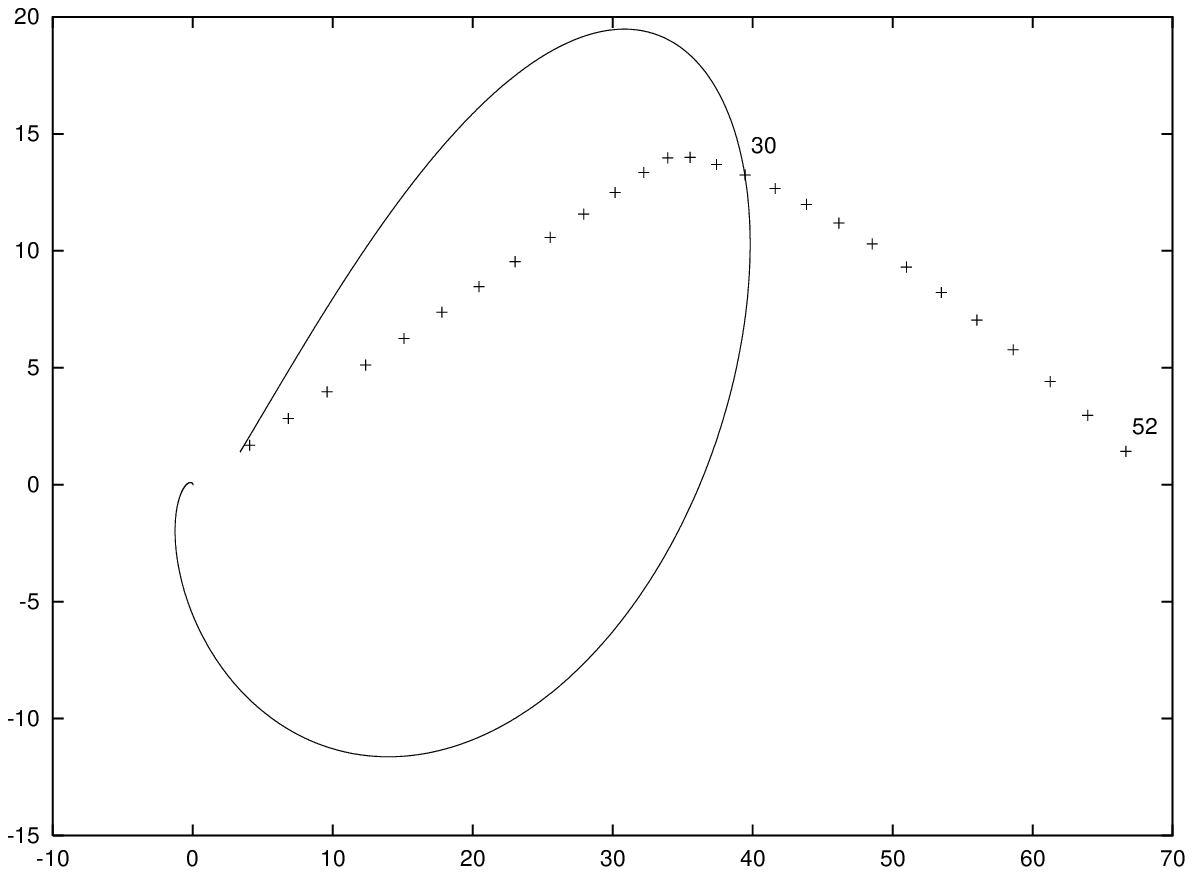}
\end{picture}
\end{figure}

\begin{center}
Figure 2. Potential and \EVs \ of $H_c$ (left) and $\tilde{H}_c$
(right)
\end{center}

The graphs of the true \EFs \ also suggest that for
a certain range of $m$ ($28 \leq m \leq 40$ and $30 \leq m \leq 44$
for $H_c$ and $\tilde{H}_c$ respectively)
the solution of (\ref{evp1}) has to be sought as a linear combination
of two JWKB modes rather than a single function of type (\ref{appr}).
These modes are determined by the corresponding values  of $(a,\,\y)$.
Of course the problem (\ref{evp1}) itself can have only
one $\cL$ eigenfunction for each \EV, but the peculiarities of the solutions
indicate that the existence of more than one JWKB approximation
has great significance for our problem.

\subsection{Cut-Off Technique}

There is a possibility that
$\Re(-i\y \xi_0(s)^{-2})>0$ for some range of $s$
for particular JWKB modes. Therefore, we cannot
use the formula (\ref{appr}) to approximate $\cL$-functions
$f$ on the whole $\R_+$. In fact, most of the modes defined by
(\ref{appr}) in our examples, although good around their centres
where $s$ is small enough, either increase as $s\rinf$ or change
rapidly near $s=-a$.

When this problem occurs, we modify the original
functions $y$ defined by (\ref{appr}) putting
\begin{equation}
\frac{y_1'(x)}{y_1(x)} =
\frac{y_1'(s_1)}{y_1(s_1)}, \ s \geq s_1,
\qquad \frac{y_2'(x)}{y_2(x)} =
\frac{y_2'(s_2)}{y_2(s_2)}, \ s \leq s_2,
\label{cut}
\end{equation}
where the points $s_1,\,s_2$ are chosen properly for each JWKB function.
Otherwise the functions are not changed anywhere.
Our aim is to replace the formerly defined modes by decaying
exponentials beyond $[s_2,s_1]$ if necessary; we shall keep
the same notation $y(x)$ for the modified functions in what follows.

It is clear that our results depend on the choice of these $s_k$.
Generally, in practical computations we have to vary two parameters
to minimise the approximation error for each particular eigenfunction
(see the next subsection for corresponding formulae).
To get optimal results one can take the initial values of $s_1$
and $s_2$ reasonably close to the corresponding local minima
of $|y_1|$ and $|y_2|$ respectively, replacing $y_k$ by
exponentially decaying functions  according to (\ref{cut}).
Then the values of $s_k$ providing the best possible
approximations among those given by the above formulae
are found by a standard minimisation procedure.

However, this general routine can be substituted by a simple
although not quite universal method.
We observe that in the examples considered here a
possible way of choosing $s_k$ is as follows. For the original
functions $y_k$ we evaluate
$$
\frac{y'_k(x)}{y_k(x)}\,=\,- i\y_k\xi_k(x) -
\frac{\xi'_k(x)}{2\xi_k(x)},
$$
then take
$$
s_k={\rm argmax} \left|\Re \frac{y'_k(x)}{y_k(x)}\right|.
$$
The reasoning leading to this choice, although not rigorous, may
prove to work in a generic situation.
Replacing $y'_k/y_k$ by constants beyond the points
$s_k$ as in (\ref{cut}), we obtain
our final modes which are almost orthogonal to one another.
Their suitable linear combinations approximate the actual \EFs \
of our operators fairly well. This is the subject of the remaining
subsections.

\subsection{Constructing Eigenfunctions}

In a generic situation when equation (\ref{ei}) has
several real roots of proper signs,
we have to approximate the relevant \EF \ by
the linear combination of the corresponding JWKB functions:
\begin{equation}
\phi(x) = \sum_{k=1}^n c_k y_k(x).
\label{comb}
\end{equation}
Here $y_k(x)$ are basic JWKB modes related to
$(a_k,\y_k)$. The number of modes entering (\ref{comb})
can be arbitrary; in our examples $n=1$ or 2.
The coefficients $c_k$ are to be determined.

One approach to computing these constants is as follows. Let
$c_k$ minimise the value of
\begin{equation}
\Delta:=\frac{\| f(x) - \phi(x)\|}{\|f\|}.
\label{de}
\end{equation}
Denote $d(x):=f(x) - \phi(x)$; a direct calculation gives
us
\begin{equation}
\|d\|^2=\|f(x) - \sum_{k=1}^{n} c_k y_k(x)\|^2 = \|f\|^2 - 2
\Re\left(\sum_{k=1}^n
c_k \bar{u}_k + \sum_{j,k=1}^n a_{jk} c_k \bar{c}_j \right)
\label{z}
\end{equation}
where
$$
u_k = \inn{f \bar{y}_k} \ , \qquad a_{jk} = \inn{y_k \bar{y}_j} \ .
$$
Differentiating (\ref{z}) with respect to $c_k$ and putting the
variation equal to zero, we obtain the desired constants:
\begin{equation}
C = A^{-1} U, \ \ C=(c_1,\ldots,c_n)^T, \ \
U=(u_1,\ldots,u_n)^T, \ \ A=(a_{jk}),
\ j,k=1,\ldots,n.
\label{c}
\end{equation}

The integrals involved into (\ref{c}) are easily computed by means of
the routine described in \cite{AD} in terms of $f'/f$, $y'_k/y_k$.
This procedure does not require computing the \EFs \
themselves: it is only the values of
$u_k,\,a_{jk},\,j,k=1,\ldots,n$,
that we have to calculate. Then the constants $c_k,\,k=1,\ldots,n$, are found from (\ref{c}).
Obviously, for a given \EF \ our $c_1,\,c_2$ must depend on the normalisation of
the functions $f$ and $y_k$ while the value of $|c_2/c_1|$ must not.
We take $\|f\|=\|e_1\|=\|e_2\|=1$
to make the matrix $A$ well-conditioned and then repeat computations for several
different values of  $\|e_k\|$ to double-check our results (we leave
$\|f\|=1$ throughout the paper).
The procedure developed in \cite{AD} turns out to be
suitable for our purposes and provides reliable answers.

Let us tabulate $c:=|c_2/c_1|$ for the examined \EVs \
of the operators $H_c$ and $\tilde{H}_c$
(as above, $m$ denotes the number of an \EV).

\vspace{.1in}

Table 5. Values of $c$ for the operators $H_c$ and $\tilde{H}_c$

\vspace{.1in}
\begin{tabular}{|c|c|c|c|c|c|c|} \hline
   $m$&  $H_c$ & $\tilde{H}_c$ &
&  $m$&  $H_c$ & $\tilde{H}_c$ \\ \hline
28 &  0.0447    & & \ & 42 &  277.2584  & 13.5057\\ \hline
30 &  0.0928    &  0.0328 & \ & 44 &  1434.62   & 48.2824 \\ \hline
32 &  0.2083    &  0.0621 & \ & 46 &  7615.34   & 177.1347 \\ \hline
34 &  0.7876    &  0.1420 & \ & 48 & 41623.37   & 660.8588  \\ \hline
36 &  2.9866    &  0.3890 & \ & 50 & 354757.06  & 2486.87  \\ \hline
38 &  13.4820   & 1.1897  & \ & 52 &            & 10590.02 \\ \hline
40 &  65.6680   & 3.9086  & & \ & & \\ \hline
\end{tabular}

\vspace{.1in}

The contents of the above table is illustrated
qualitatively by the plots
of the \EFs \ (see figure~3). Indeed, for a certain range of
\EVs \ the relevant
\EFs \ have two distinct components clearly  seen in the plots.
For lower \EVs \ the first mode $y_1$ is dominating, while
the contribution of $y_2$ is more significant for
the high energy spectrum. All our numerical results
including the plots of the \EFs \ indicate the transition from
the first leading JWKB function
$y_1$ to $y_2$ in formula (\ref{comb}). This phenomenon
is observed for both examples studied here.

In the previous subsection we mentioned the parameters
$s_k$ entering the modified functions $y_k$. The
formulae of this subsection contain these parameters
fixed at some starting values. After minimising
$\Delta$ with respect to $c_k$ as described we minimise the obtained
function with respect to $s_k$, and thus
get the final results to be tabulated in the next subsection.

\subsection{The Efficiency of Approximations}

Having found $c_k,\,k=1,2$, we then compute
$\Delta$ defined by (\ref{de}) in terms of the calculated coefficients
$u_k,\,a_{jk}$. The values of $\Delta$ tabulated below characterise the
accuracy of the JWKB approximations we use.

As is seen from table~6, not only do the JWKB modes
approximate \EFs \
corresponding to higher \EVs \ but they also turn out
to be surprisingly efficient even for the lower part of
the spectrum. The plots of $|d_m|$ in comparison with $|f_m|$
(see figure~4) also indicate the
effectiveness of the approximations obtained by our approach.
The corrections $d_m$ for other values of $m$ either look similar to those
shown in figure~4 or are invisible compared to the corresponding \EFs.

These numerical results complete our study of
eigenfunction approximations.

\vspace{.1in}
Table~6. Minimal values of $\Delta$

\vspace{.1in}
\begin{tabular}{|c|c|c|c|c|c|c|} \hline
  $m$ & $H_c$ & $\tilde{H}_c$ & \ &
  $m$ & $H_c$ & $\tilde{H}_c$ \\ \hline
10 &  0.132627 & 0.153114
& \ &32 & 0.015264& 0.015016\\ \hline
12 &   0.042214& 0.048187
& \ &34 & 0.017776& 0.015332\\ \hline
14 &   0.029275& 0.032804
& \ &36 & 0.014177& 0.020596 \\ \hline
16 &   0.021119& 0.024966
& \ &38 & 0.009434& 0.025885\\ \hline
18 &   0.015969& 0.019302
& \ &40 & 0.006557& 0.025631\\ \hline
20 &   0.012366& 0.011576
& \ &42 & 0.004796& 0.012372\\ \hline
22 &   0.011045& 0.009765
& \ &44 & 0.003873& 0.007280\\ \hline
24 &   0.010724& 0.009110
& \ &46 & 0.002828& 0.005568\\ \hline
26 &   0.014866& 0.008487
& \ &48 & 0.002236& 0.004359\\ \hline
28 & 0.015133&  0.009798
& \ &50 & 0.001732& 0.003464\\ \hline
30 & 0.015232& 0.012214
& \ &52 &         & 0.003133\\ \hline
\end{tabular}
\vspace{.1in}

\pagebreak
\begin{figure}[h]
\vspace{1cm}
\begin{picture}(200,400)(0,0)
\includegraphics{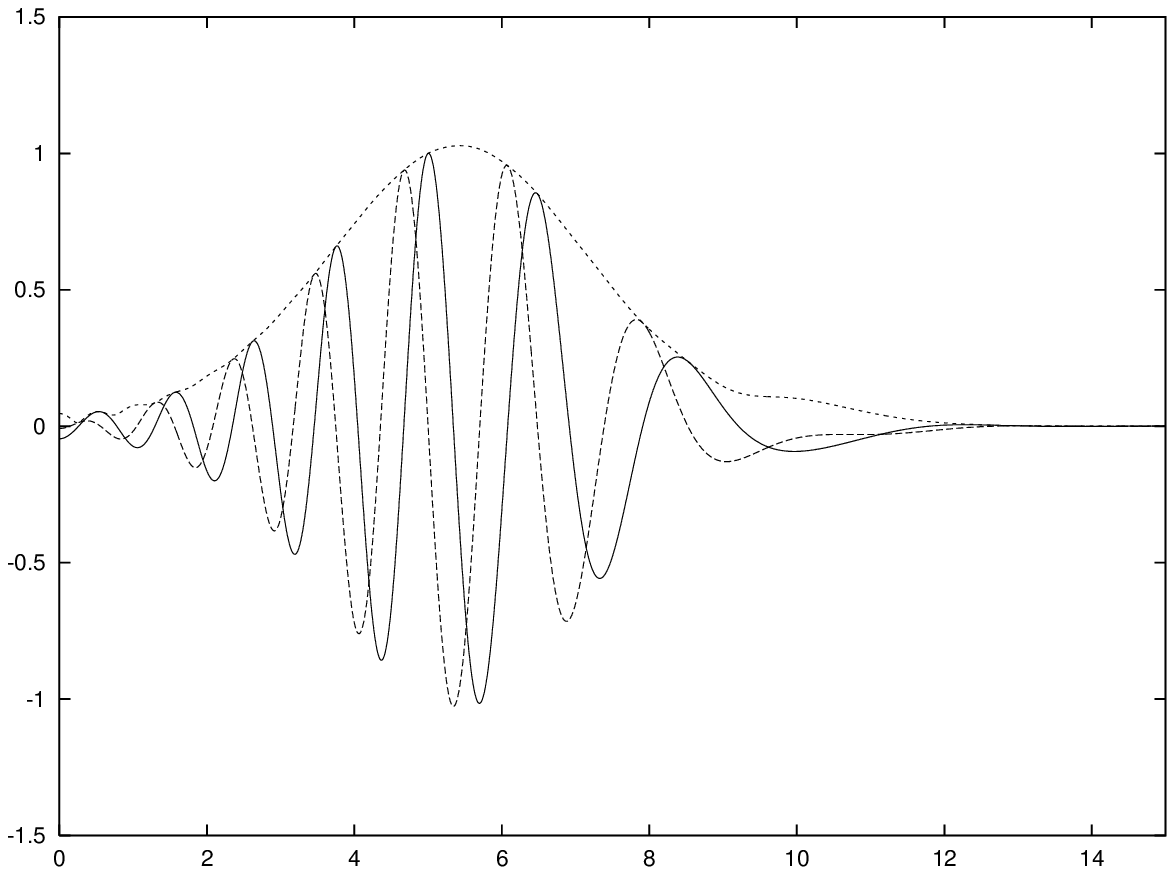}
\includegraphics{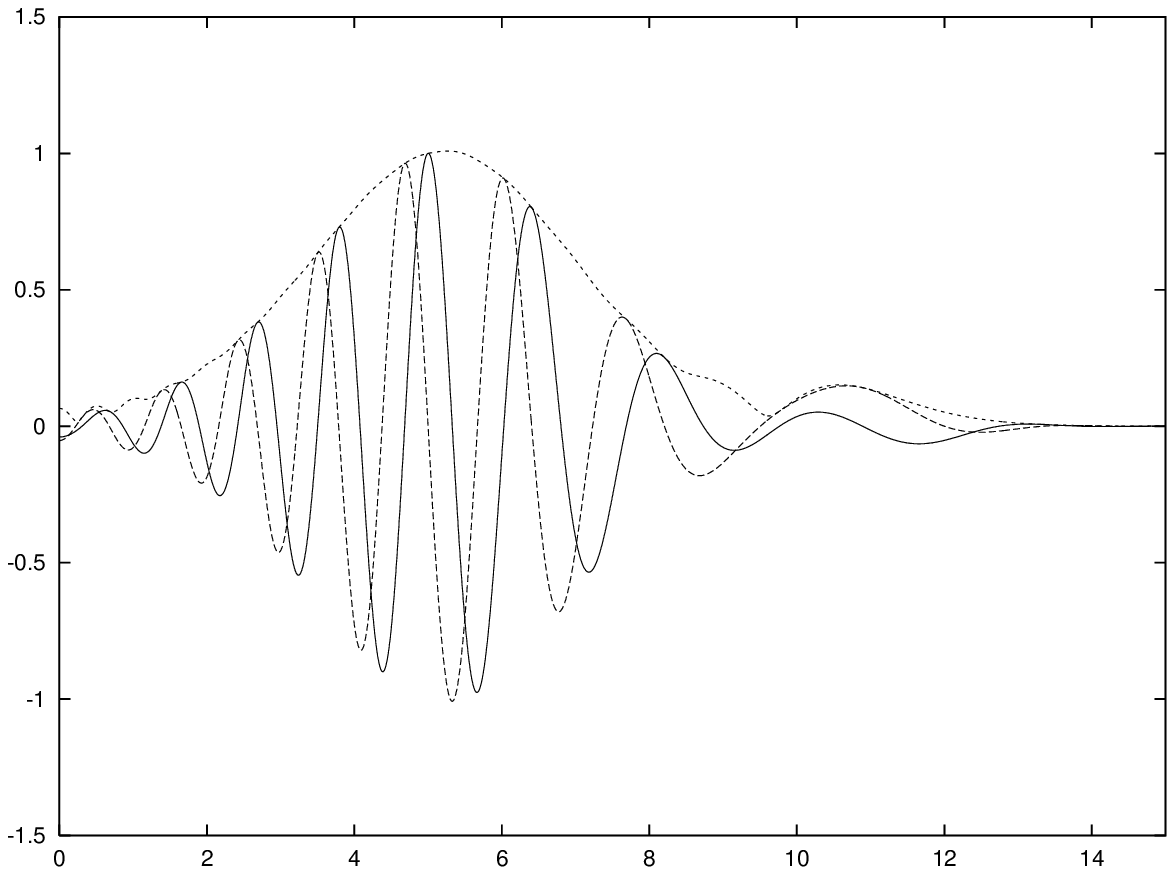}
\includegraphics{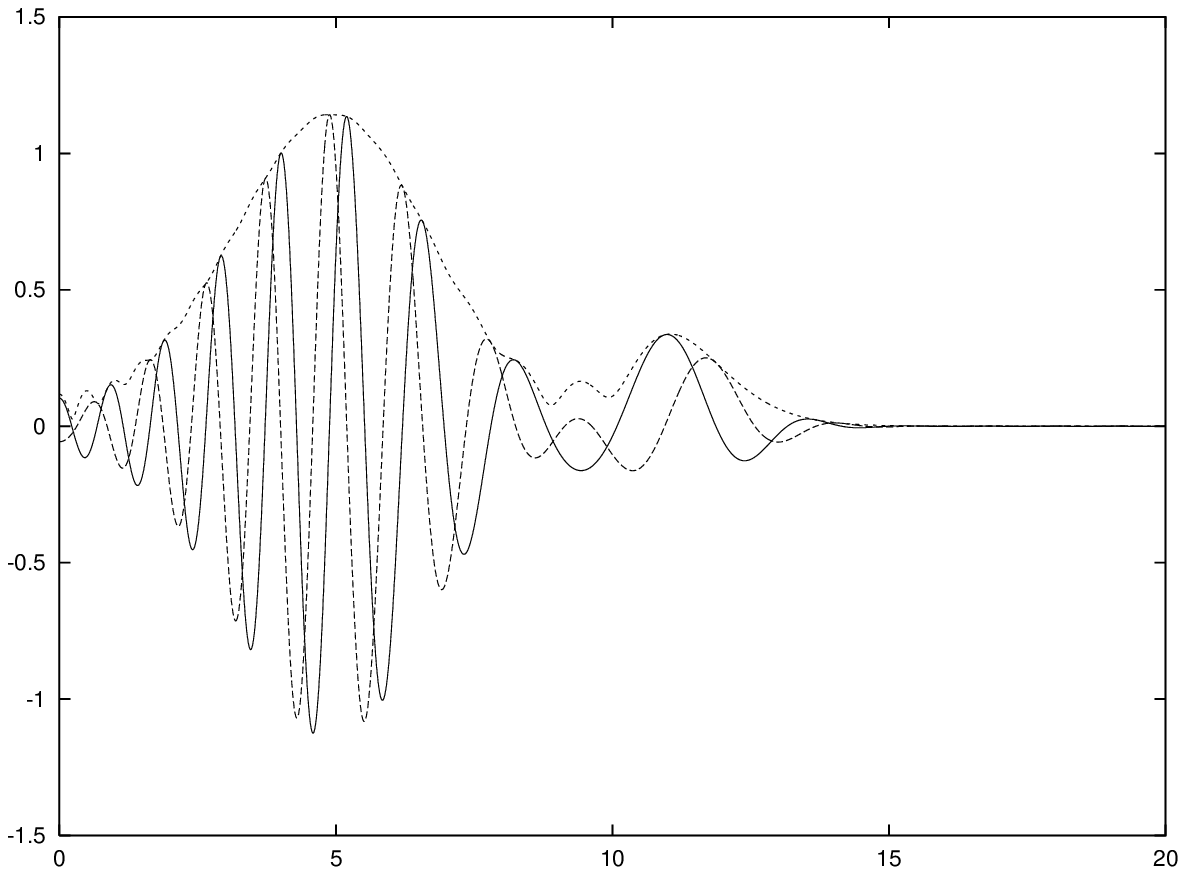}
\includegraphics{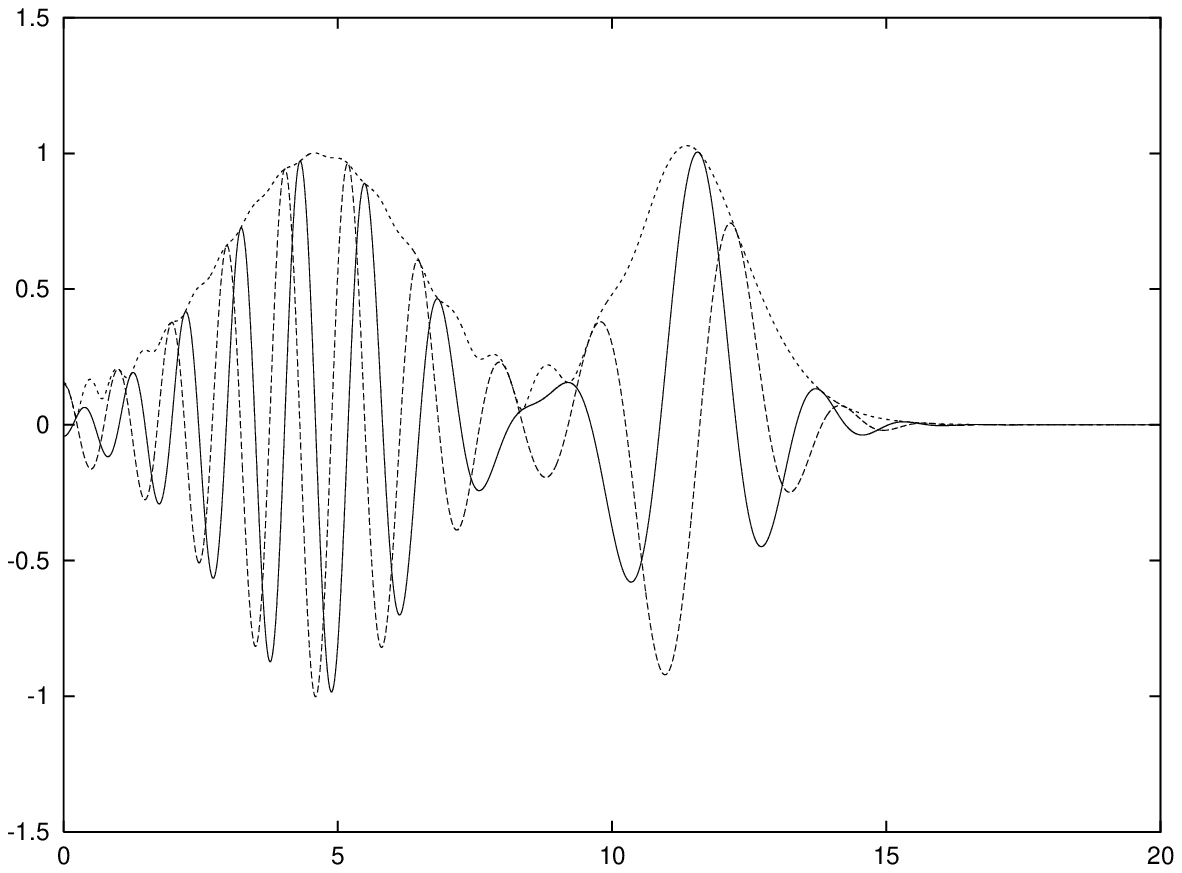}
\includegraphics{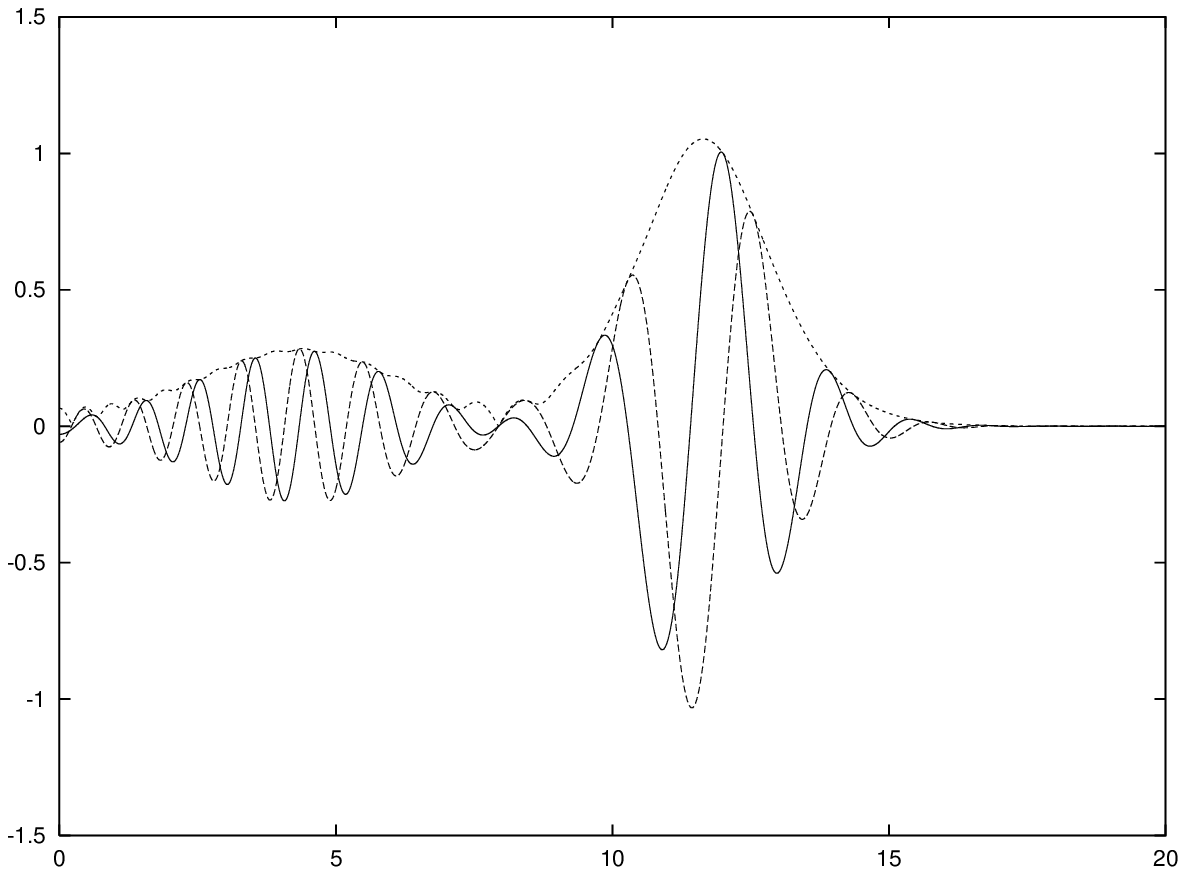}
\includegraphics{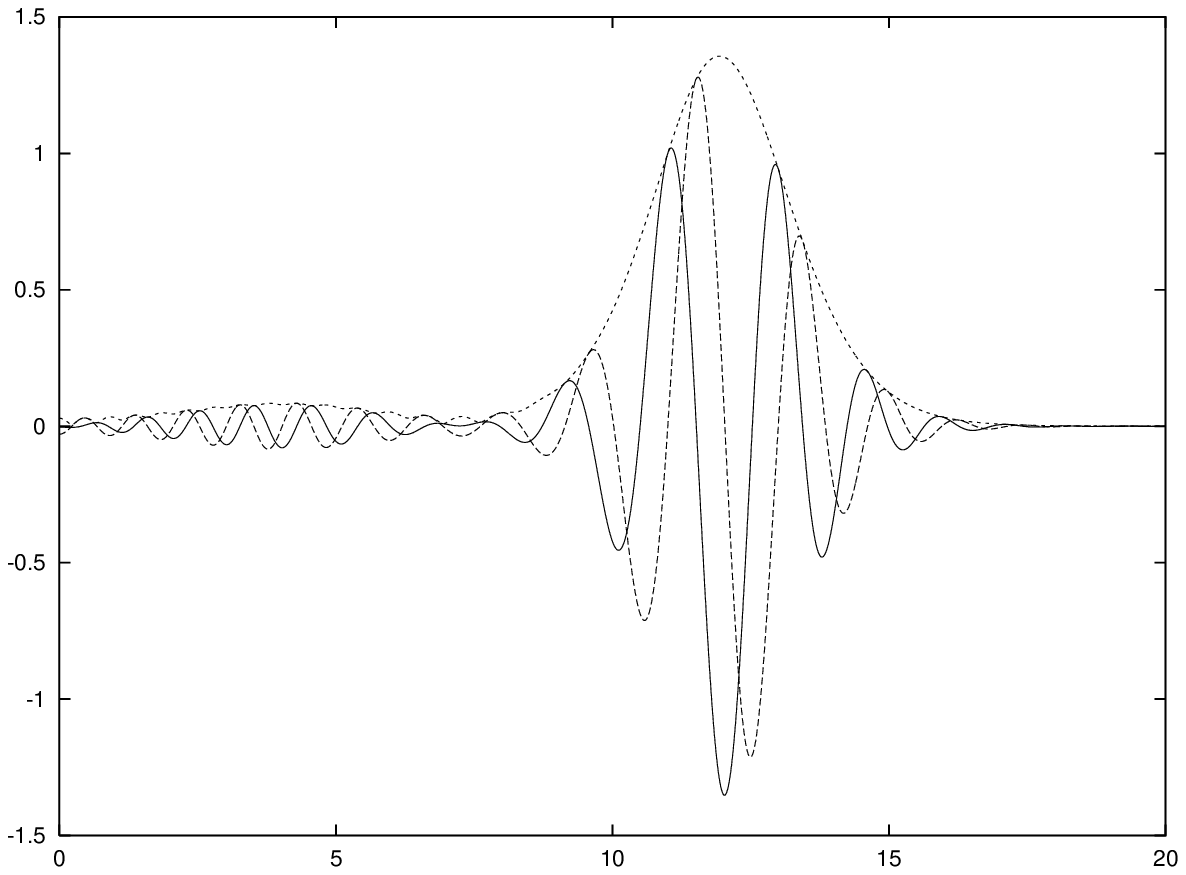}
\end{picture}
\end{figure}

\vspace{0.3in}
\begin{center}
Figure 3. Eigenfunctions \ of $H_c$: $\Re{f_m},\,\Im{f_m},\,|f_m|$. From left to right, top to bottom: $m=28,\,30,\,32,\,34,\,36,\,38$
\end{center}

\pagebreak
\begin{figure}[h]
\vspace{1cm}
\begin{picture}(200,100)(0,0)
\includegraphics{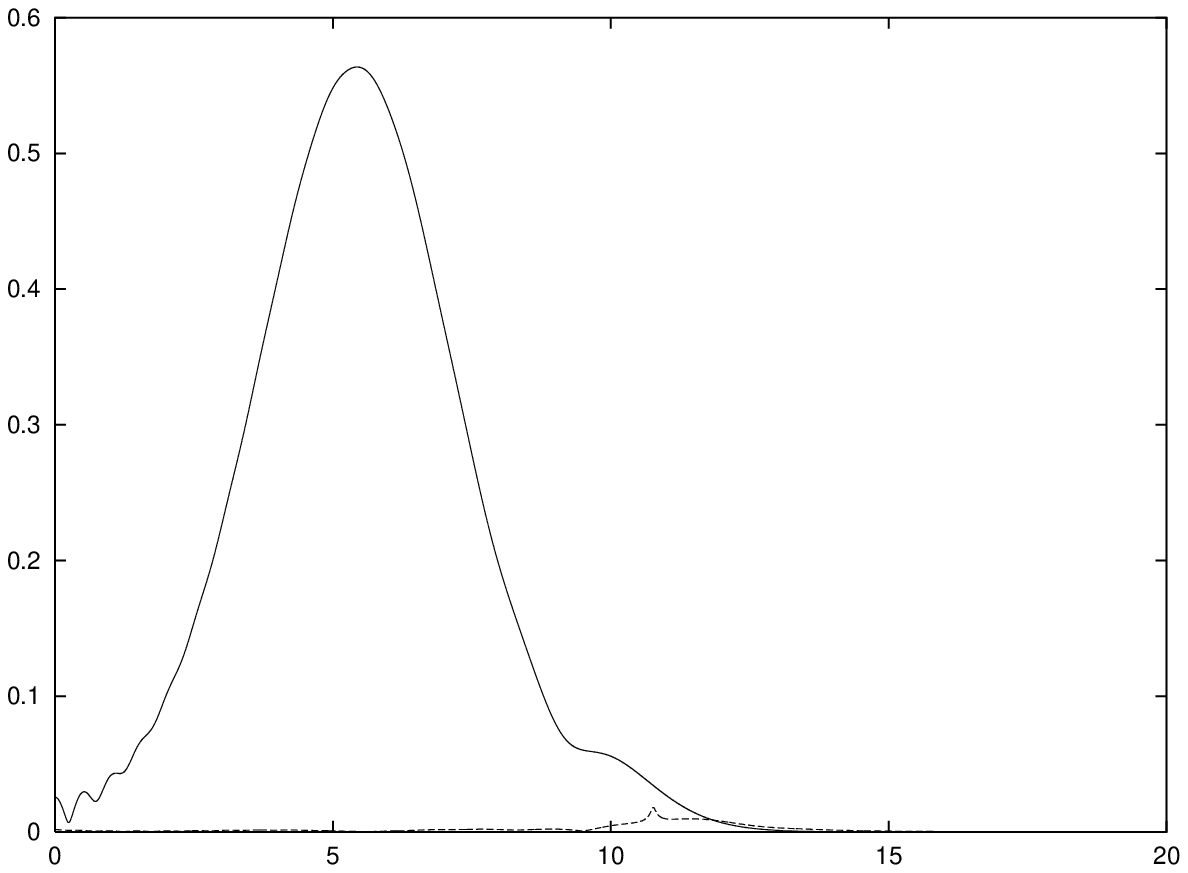}
\includegraphics{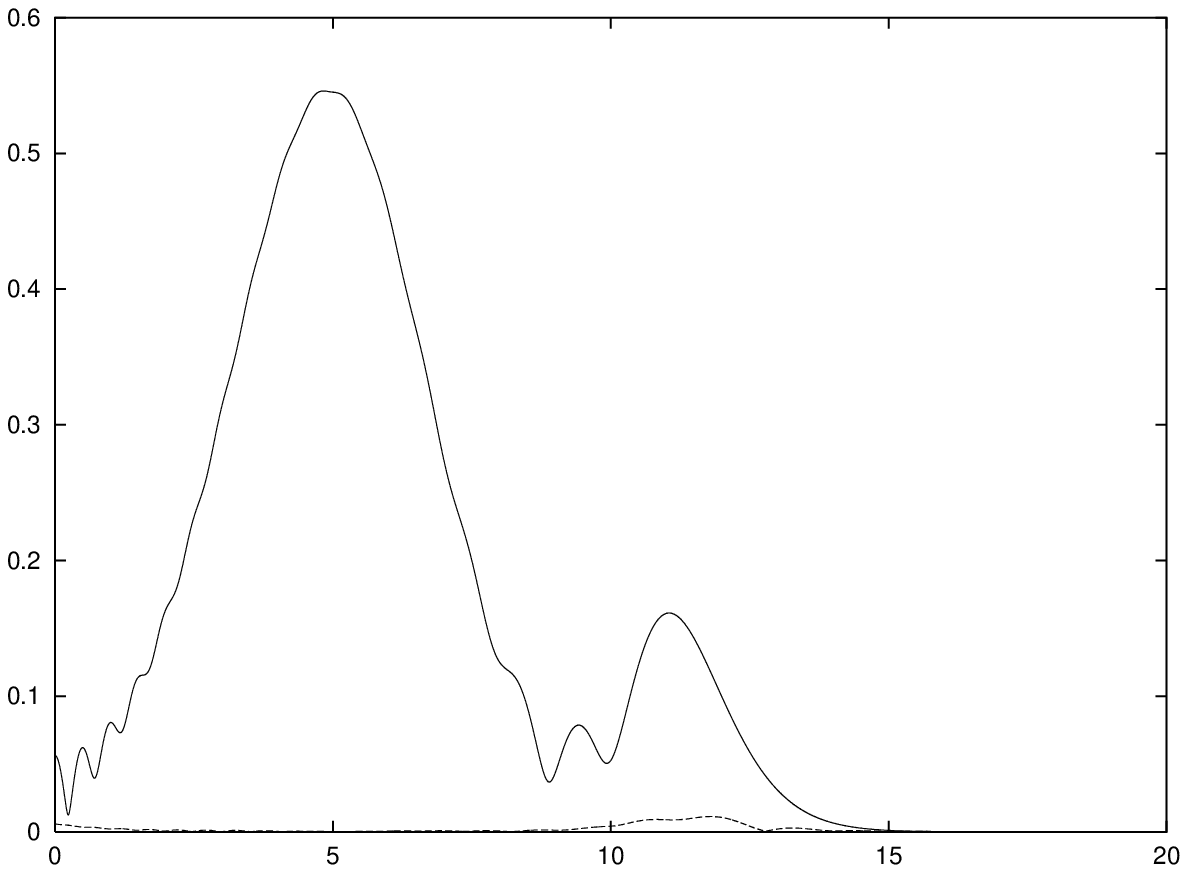}
\end{picture}
\vspace{1in}
\begin{center}
Figure 4. Plots of $|f_m|$ and $|d_m|$ for $H_c$. From left to right: $m=28,\,32$
\end{center}
\end{figure}

\pagebreak

\end{document}